\documentclass[smallextended]{svjour3}   
\smartqed  
\usepackage{graphicx}
\usepackage{amsmath,amssymb,latexsym,enumerate,mathrsfs}
\usepackage{hyperref}

\DeclareMathOperator{\diag}{diag}
\newtheorem{assumption}[theorem]{\bf Assumption}
\def\qed{\hfill $\Box$}
\def\etal{\mbox{et al.}}
\usepackage{graphicx}
\journalname{JOTA}

\begin{document}

\title{A further study on the opioid epidemic dynamical model with random perturbation}

\titlerunning{A further study on the opioid epidemic dynamical model} 

\author{Getachew K. Befekadu \and
Quanyan Zhu
}

\institute{Getachew K. Befekadu \at
              Department of Mechanical and Aerospace Engineering\\
              University of Florida - REEF, 1350 N. Poquito Rd, Shalimar, FL 32579, USA. \\
              Tel.: +1 850 833 9350\\
              Fax: +1 850 833 9366\\
              \email{gbefekadu@ufl.edu}    
               \and
           Quanyan Zhu \at
           Department of Electrical and Computer Engineering\\
           Tandon School of Engineering\\
           New York University\\
           5 MetroTech Center, LC 200A, Brooklyn, NY 11201, USA.\\
            \email{quanyan.zhu@nyu.edu}
}

\date{Received: date / Accepted: date}

\maketitle

\begin{abstract}
In this paper, we consider an opioid epidemic dynamical model with random perturbation that typically describes the interplay between regular prescription use, addictive use, and the process of rehabilitation from addiction and vice-versa. In particular, we provide two-sided bounds on the solution of the transition density function for the Fokker-Planck equation that corresponds to the opioid epidemic dynamical model, when a random perturbation enters only through the dynamics of the susceptible group in the compartmental model. Here, the proof for such bounds basically relies on the interpretation of the solution for the transition density function as the value function of a certain optimal stochastic control problem. Finally, as a possible interesting development in this direction, we also provide an estimate for the attainable exit probability with which the solution for the randomly perturbed opioid epidemic dynamical model exits from a given bounded open domain during a certain time interval. Note that such qualitative information on the first exit-time as well as two-sided bounds on the transition density function are useful for developing effective and fact-informed intervention strategies that primarily aim at curbing opioid epidemics or assisting in interpreting outcome results from opioid-related policies.
\end{abstract}

 

\section{Introduction} \label{S1}

In recent years, the Internet and its continued promotion of prescription opioid abuse, as well as the inappropriate physicians prescribing practices have exacerbated the opioid epidemic by making opioids more accessible or increased the supply of unused opioids available for further misuse. In order to address this modern era epidemic plague -- mainly driven by opioid addiction, which degrades health, devastates families and reduces productivity at a huge societal and economic cost -- a number of federal and state agencies throughout the United States have implemented a wide range of opioid-related policies\footnote{Including the Ryan Haight Online Pharmacy Consumer Protection Act of 2008 which prohibited the Internet distribution of controlled substances without a valid prescription \cite{r5}; see also \cite{DowHC16} for CDC guideline for prescribing opioids for chronic pain--United States, 2016.} that are primarily aimed at curbing prescription opioid abuse, establishing guidelines to prevent inappropriate prescribing practices, developing abuse deterrents or preventing drug diversion mechanisms \cite{r1}, \cite{r2},  \cite{r3} and \cite{r4}. On the other hand, only a few studies have been reported on the need for effective and fact-informed intervention strategies, based on mathematical theory of epidemiology for infectious diseases, with the intent of better understanding the dynamics of the current serious opioid epidemic (e.g., see \cite{r6} and \cite{r7} in context of exploring the dynamics of drug abuse epidemics, focusing on the interplay between the different opioid user groups and the process of rehabilitation and treatment from addiction; see \cite{r8}, \cite{r9}, \cite{r10}, \cite{r11} or \cite{r12} for additional studies, but in the context of heroin epidemics that resembling the classic susceptible-infected-recovered (SIR) model, based on the work of \cite{r13}; see also \cite{CasJR97} for an interesting discussion of mathematical models for the dynamics of tobacco use, recovery and relapse).

In this paper, without attempting to give a literature review, we consider an opioid epidemic dynamical model that describes the interplay between regular prescription use, addictive use, and the process of rehabilitation from addiction and vice-versa (e.g., see \cite{r6} for additional discussions). In particular, we are mainly interested in understanding how a random noise propagates through the opioid epidemic dynamical model, when the random perturbation enters only through the dynamics of the susceptible group in the compartmental model. Here, we provide two-sided bounds on the solution of the transition density function for the corresponding Fokker-Planck equation that is associated with randomly perturbed opioid epidemic dynamical model. The proof for such bounds basically relies on the interpretation of the solution for the density function as the value function of a certain optimal stochastic control problem. Moreover, as a possible interesting development in this direction, we also provide an estimate for the attainable exit probability with which the solution for the randomly perturbed opioid epidemic dynamical model exits from a given bounded open domain during a certain time interval.

Here, it is worth mentioning that some interesting studies how random noise may propagate through a chain of dynamical systems have been reported in literature (e.g., see \cite{BefA15a} and \cite{DelM10} mainly from mathematical point of view). The rationale behind our framework, which follows in some sense the settings of these papers, is to provide a probabilistic representation for the transition density function and further make a connection with the value function of some optimal stochastic control problem, via a logarithmic transformation (e.g., see \cite{FleSh85} or \cite{She91} for a similar argument). Note that qualitative information obtained from such a stochastic control argument are useful for developing effective and fact-informed intervention strategies that primarily aim at curbing opioid epidemics, or assisting in interpreting outcome results from opioid-related policies.\footnote{In this paper, our intent is to provide a theoretical framework, rather than providing any specific numerical simulations.}

The remainder of this paper is organized as follows. In Section~\ref{S2}, we provide some preliminary results that will be useful in the paper. In Section~\ref{S3}, we provide two-sided bounds on the solution of the density function for the Fokker-Planck equation that corresponds to the opioid epidemic dynamical model, when there is a random perturbation enters through the dynamics of the susceptible group in the compartmental model. Finally, Section~\ref{S4} provides an estimate for the attainable exit probability of the diffusion process, which is associated with randomly perturbed opioid epidemic dynamical model from a given bounded open domain during a certain time interval.

\section{Preliminaries} \label{S2}
In this section, we provide some preliminary results that will be useful later in the paper.
\subsection{Mathematical model} \label{S2(1)}
\begin{figure}[ht]
\includegraphics[width=3.15 in]{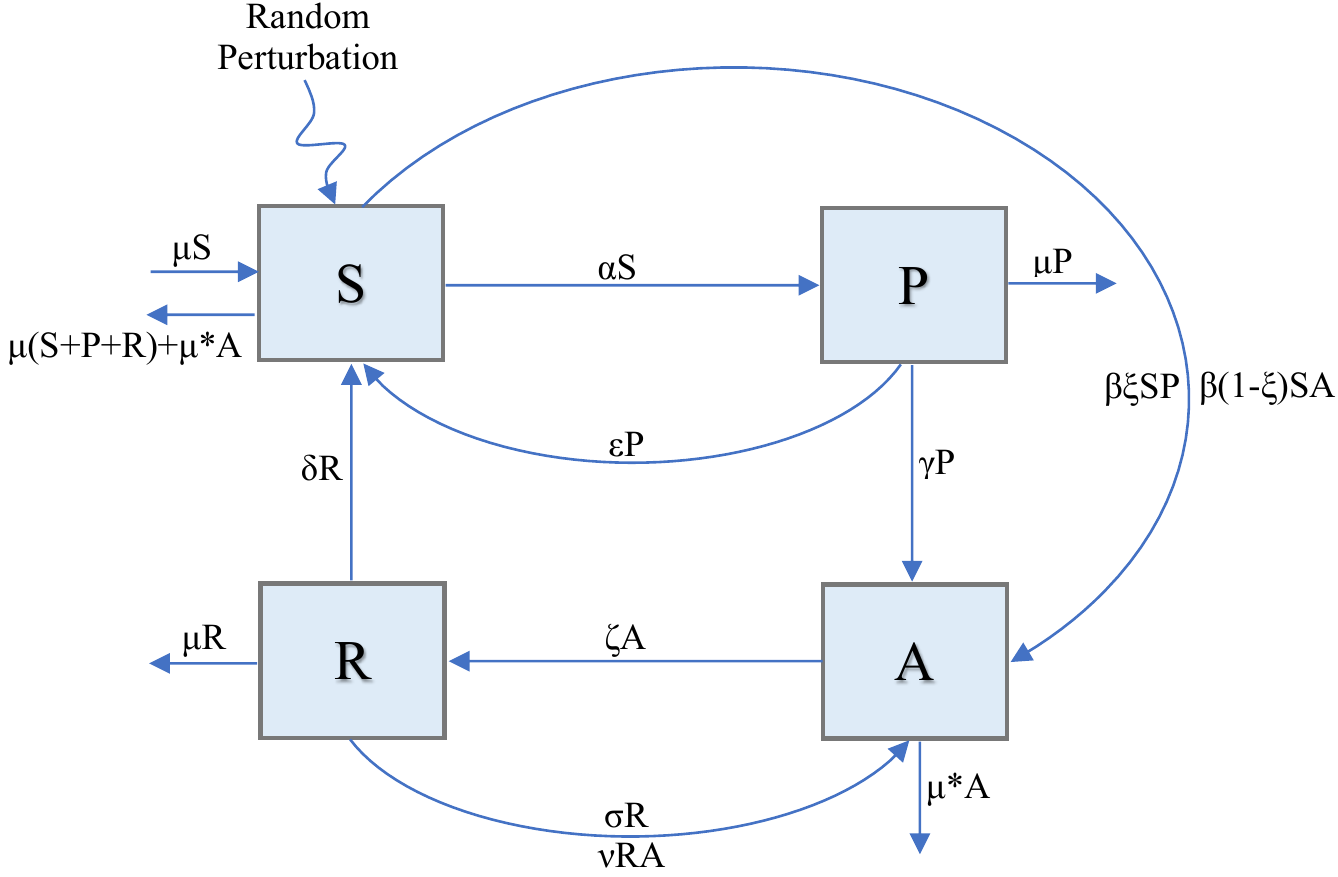}
\caption {\scriptsize A block diagram showing the relationships between the different groups in the compartmental model of opioid addiction with random perturbation (cf. Battista \etal \,\cite{r6})} \label{Fig1}
\end{figure}

In this paper, we specifically consider an opioid epidemic dynamical model that describes the interplay between regular prescription opioid use, addictive use, and the process of rehabilitation from addiction and vice-versa (e.g., see \cite{r6} for a detailed discussion). To this end, we introduce the following population groups
\begin{enumerate} [(i)]
\item {\it Susceptible group} - $S$: This group in the compartmental model includes those who are susceptible to opioid addiction, but not currently using opioids. In the compartmental model, everyone who is not in addiction treatment, already an addict, or using opioids as medically prescribed is classified as ``susceptible".
\item {\it Prescribed user group} - $P$: This group in the compartmental model is composed of individuals who have health related concerns and also have access to opioids through a proper physician's prescription, but they are not addicted to opioids. Members of this group have some inherent tendency of becoming addicted to their prescribed opioids.
\item {\it Addiction user group} - $A$: This group in the compartmental model is composed of people who are addicted to opioids. There are multiple interaction routes to this group in the compartmental model, including those routes that are bypassing the prescribed user group $P$ (see Fig~\ref{Fig1}).
\item {\it Treatment/rehabilitation} - $R$: This group in the compartmental model contains individuals who are in treatment for their addiction. Here, we include an inherent rate of falling back into addiction as well as a typical process of relapsing due to general availability of the drug. Moreover, we also assume that some of the members from the recovering group who have completed their treatment may return to being susceptible. That is, we assume that successful treatment does not imply permanent immunity to addiction (i.e., in general, an assumption based on the balance of increased risk of addiction verses increased awareness and avoidance).
\end{enumerate}
Then, using the basic remarks made above, we specify the following SIR compartmental model for the opioid epidemics which is described by the following four continuous-time differential equations 
\begin{align}
\dot{S}(t) = -\alpha S(t) &- \beta(1-\xi)S(t) A(t) - \beta \xi S(t) P(t) + \epsilon P(t) \notag \\
                        &\quad  + \delta R(t) + \mu(P(t)+R(t)) + \mu^{\ast} A(t), \label{Eq2.1a}
\end{align}
\begin{align}
\dot{P}(t) &= \alpha S(t) - (\epsilon + \gamma + \mu) P(t), \label{Eq2.1b}
\end{align}
\begin{align}
\dot{A}(t) =  \gamma P(t) + \sigma R(t) &+ \beta(1 - \xi) S(t) A(t) + \beta \xi S(t) P(t) \notag \\
                                                    & + \nu R(t) A(t) - (\zeta + \mu^{\ast}) A(t) \label{Eq2.1c}
\end{align}
and
\begin{align}
\dot{R}(t) &= \zeta A(t) - \mu R(t) A(t) - (\delta + \sigma + \mu) R(t), \label{Eq2.1d} 
\end{align}
where the normalized overall population is assumed to be constant, i.e., $1= A(t)+S(t)+R(t)+P(t)$, since the number of mortality due to opioid-related overdose is very small, when compared to the change in the total population numbers in the short term. Moreover, the followings are brief description for the parameters in the above system of equations (i.e., the system parameters in Equations~\eqref{Eq2.1a}--\eqref{Eq2.1d})
{\em
\begin{itemize}
\item $\alpha S(t)$: the rate at which people are prescribed opioids.
\item $\beta$: the total probability of becoming addicted to opioids other than by prescription.
\item $\beta(1-\xi)$: the proportion of $\beta$ caused by black market drugs or other addicts.
\item $\beta \xi$: the rate at which the non-prescribed, susceptible population begins abusing opioids due to the accessibility of extra prescription opioids, e.g., new addicts got the drug from a friend or relative's prescription.
\item $\epsilon$: the rate at which people come back to the susceptible group after being prescribed opioids.
\item $\delta$: the rate at which people come back to the susceptible group after successfully finishing treatment. Despite having completed rehabilitation, we assume people are susceptible to addiction for life.
\item $\mu$: the natural death rate.
\item $\mu^{\ast}$: the (enhanced) death rate for addicts ($\mu$ plus overdose rate).
\item $\gamma$: the rate at which the prescribed opioid users fall into addiction.
\item $\zeta$: the rate at which addicted/dependent opioid users enter the treatment/ rehabilitation process.
\item $\nu$: the rate at which users during the treatment fall back into addictive drug use due to the availability of prescribed painkillers from friends or relatives.
\end{itemize}}

Note that the normalized overall population is assumed constant (which is set to unity). Then, with $P(t) = 1 - S(t) - A(t) - R(t)$, we can reduce the above system of equations in Equations~\eqref{Eq2.1a}--\eqref{Eq2.1d} as follows
\begin{eqnarray}
\left.\begin{array}{l}
\dot{S}(t) = -\alpha S(t) -\beta(1-\xi)S(t) A(t) - \beta \xi S(t) (1 - S(t) - A(t) - R(t))  \\
\hspace{0.9in} + (\epsilon + \mu)(1 - S(t) - A(t) - R(t)) + (\delta + \mu) R(t) + \mu^{\ast} A(t)\\
\dot{A}(t) =  \gamma (1 - S(t) - A(t) - R(t)) + \sigma R(t) + \beta(1 - \xi) S(t) A(t) \\
 \hspace{0.9in} + \beta \xi S(t) (1 - S(t) - A(t) - R(t))  + \nu R A - (\zeta + \mu^{\ast})A(t) \\
\dot{R}(t) = \zeta A(t) - \mu R(t) A(t) - (\delta + \sigma + \mu) R(t)
\end{array}\right\}. \label{Eq2.2}
\end{eqnarray}
In order to facilitate our presentation, we adopt the following change of variables: $S \rightarrow x_1$, $A \rightarrow x_2$ and $R \rightarrow x_3$. Then, with minor abuse of notation, the system of equations in Equation~\eqref{Eq2.2} can be further rewritten as follows
\begin{eqnarray}
\left.\begin{array}{l}
\dot{x}_1(t) = f_1(t, x_1, x_2, x_3)\\
\dot{x}_2(t) = f_2(t, x_1, x_2, x_3) \\
\dot{x}_3(t) = f_3(t, x_2, x_3)
\end{array}\right\} \label{Eq2.2b}
\end{eqnarray}
where the functions $f_1$, $f_2$ and $f_3$ are given by
\begin{align*}
& f_1(t, x_1, x_2, x_3) \\
 &\quad \quad  =  -\alpha x_1(t) -\beta(1-\xi)x_1(t) x_2(t) - \beta \xi x_1(t) (1 - x_1(t) - x_2(t) - x_3(t)) \\
                                &\quad \quad\quad  + (\epsilon + \mu)(1 - x_1(t) - x_2(t) - x_3(t)) + (\delta + \mu) x_3(t) + \mu^{\ast} x_2(t),
\end{align*}
\begin{align*}
f_2(t, x_1, x_2, x_3) =& \gamma (1 - x_1(t) - x_2(t) - x_3(t)) + \sigma x_3(t) + \beta(1 - \xi) x_1(t) x_2(t)  \\
                                 & \quad + \beta \xi x_1(t) (1 - x_1(t) - x_2(t) - x_3(t)) + \nu x_3(t) x_2(t) \\
                                 & \quad \quad \quad - (\zeta + \mu^{\ast}) x_2(t)
\end{align*}
and
\begin{align*}
f_3(t, x_2, x_3) = \zeta x_2(t) - \mu x_3(t) x_2(t) - (\delta + \sigma + \mu) x_3(t).
\end{align*}
respectively.

In what follows, we assume that a random noise enters only through the dynamics of the susceptible group in Equation~\eqref{Eq2.2} and is then subsequently propagated to the other groups in the compartmental model (see also in Fig~\ref{Fig1}). To this end, we consider the corresponding system of stochastic differential equations (SDEs), i.e.,
\begin{eqnarray}
\left.\begin{array}{l}
dX_1(t) = f_1(t, X_1(t), X_2(t), X_3(t)) dt + \hat{\sigma}(t, X_1(t), X_2(t), X_3(t)) dW(t)\\
dX_2(t) = f_2(t, X_1(t), X_2(t), X_3(t)) dt \\
dX_3(t) = f_3(t, X_2(t), X_3(t)) dt
\end{array}\right\} \label{Eq2.6}
\end{eqnarray}
where $\bigl(W(t)\bigr)_{t \ge 0}$ is a one-dimensional Brownian motion, $\bigl(X_1(t), X_2(t), X_3(t)\bigr)_{t \ge 0}$ being an $\mathbb{R}^3$-valued degenerate diffusion process, and $\hat{\sigma}$ and $\hat{\sigma}^{-1}$ are assumed to be bounded functions. Moreover, if we denote by a bold letter a quantity in $\mathbb{R}^3$, for example, the solution in Equation~\eqref{Eq2.6} is denoted by $\bigl(\mathbf{X}(t)\bigr)_{t \ge 0} = \bigl(X_1(t), X_2(t), X_3(t)\bigr)_{t\ge 0}$, then we can rewrite Equation~\eqref{Eq2.6} as follows
\begin{align}
d \mathbf{X}(t) = \mathbf{F} (t, \mathbf{X}(t)) dt + B \hat{\sigma}(t, \mathbf{X}(t)) dW(t), \label{Eq2.7}
\end{align}
where $\mathbf{F} = \bigl[f_1, f_2, f_3\bigr]^T$ is an $\mathbb{R}^3$-valued function and $B$ stands for a column vector that embeds $\mathbb{R}$ into $\mathbb{R}^3$, i.e., $B = [1, 0, 0]^T$. 

Note that the backward operator for the diffusion process $\mathbf{X}(t)$, when applied to a certain function $\upsilon(t,\mathbf{x})$, is given by
\begin{align}
\partial_{t} \upsilon(t,\mathbf{x}) + \mathcal{L} \upsilon(t,\mathbf{x}) = \partial_{t} \upsilon(t,\mathbf{x}) &+  \frac{1}{2} \operatorname{tr}\Bigl \{a(t,\mathbf{x}) D_{x_1}^2 \upsilon(t,\mathbf{x}) \Bigr\}  \notag \\
& + \sum\nolimits_{i=1}^3 f_{i}(t,\mathbf{x}) D_{x_i} \upsilon(t,\mathbf{x}),  \label{Eq2.8}
\end{align}
where $a(t,\mathbf{x})=\hat{\sigma}(t,\mathbf{x})\,\hat{\sigma}^T(t,\mathbf{x})$, $D_{x_i}$ and $D_{x_1}^2$ (with $D_{x_1}^2 = \bigl({\partial^2 }/{\partial x_1 \partial x_1} \bigr)$) are the gradient and the Hessian (w.r.t. the variable $x_i$, for $i \in \{1,2,3\}$), respectively.

Let $D \subset \mathbb{R}^3$ be a given bounded open domain, with smooth boundary $\partial D$ (i.e., $\partial D$ is a manifold of class $C^2$). With $\Omega = (0, T) \times D$, let us denote by $C^{\infty}(\Omega)$ the spaces of infinitely differentiable functions on $\Omega$, and by $C_0^{\infty}(\Omega)$ the space of the functions $\phi \in C^{\infty}(\Omega)$ with compact support in $\Omega$. A locally square integrable function $\upsilon(s, \mathbf{x})$ on $\Omega$ is said to be a distribution solution to the following equation
\begin{eqnarray}
\partial_{t} \upsilon(t,\mathbf{x}) + \mathcal{L} \upsilon(t,\mathbf{x}) = 0,  \label{Eq2.9}
\end{eqnarray}
if, for any test function $\phi \in C_0^{\infty}(\Omega)$, the following holds true
\begin{eqnarray}
 \int_{\Omega} \Bigl(-\partial_{t} \phi(t,\mathbf{x}) + \mathcal{L}^{\ast} \phi(t,\mathbf{x})  \Bigr)\upsilon(t,\mathbf{x}) d \Omega = 0, \label{Eq2.10}
\end{eqnarray}
where $d \Omega$ denotes the Lebesgue measure on $\mathbb{R} \times \mathbb{R}^3$ and $\mathcal{L}^{\ast}$ is an adjoint operator corresponding to $\mathcal{L}$.

The following statements are standing assumptions that hold throughout the paper.
\begin{assumption} \label{AS1} ~\\\vspace{-3mm}
\begin{enumerate} [(a)]
\item The functions $\hat{\sigma}(s,\mathbf{x})$ and $\hat{\sigma}^{-1}(s,\mathbf{x})$ are bounded $C^{\infty} \bigl((0, \infty) \times \mathbb{R}^3\bigr)$-functions, with bounded first derivatives. Moreover, the least eigenvalue of $a(t,\mathbf{x})$ is uniformly bounded away from zero, i.e.,
\begin{align*}
  \mathbf{y}^Ta(t,\mathbf{x}) \mathbf{y} \ge \lambda \bigl\vert \mathbf{y} \bigr\vert^2, \quad \forall \mathbf{x}, \mathbf{y} \in \mathbb{R}^{3}, \quad \forall t \ge 0,
\end{align*}
for some $\lambda > 0$.
\item The backward operator in Equation~\eqref{Eq2.8} is hypoelliptic in $C^{\infty}(\Omega_{0}^{\infty})$ (e.g., see \cite{Hor67} or \cite{Ell73}).
\end{enumerate}
\end{assumption}

\begin{remark} \label{R1}
Note that the hypoellipticity assumption implies that the degenerate diffusion process $\mathbf{X}(t)$ has a transition probability density with a strong Feller property. Here, we also assume that the system of equations in Equation~\eqref{Eq2.6} satisfies a weak H\"{o}rmander's condition, where the drift terms are responsible for propagating the noise through the system dynamics (e.g., see \cite{BefA15a} or \cite{DelM10}). Moreover, we observe that there is an intrinsic time-scaling property with which the system dynamics in Equation~\eqref{Eq2.6} is evolving (cf. Equations~\eqref{Eq2.1a}--\eqref{Eq2.1d}). Later in Section~\ref{S3}, we use this intrinsic property to support our claim, i.e., for proving bounds on the solution of the transition density function for the Fokker-Planck equation (cf. Equation~\eqref{Eq2.11} below) that is associated with randomly perturbed opioid epidemic dynamical model.
\end{remark}

\subsection{Connection with stochastic control problems} \label{S2(2)}
In this subsection, we first discuss the existence of a transition density function, i.e., $\bigl(p(s,t,\mathbf{x}, \mathbf{y})\bigr)_{0 \le s < t; \mathbf{x}, \mathbf{y} \in \mathbb{R}^3}$, to the solution $\bigl(\mathbf{X}(t)\bigr)_{t \ge 0}$ of Equation~\eqref{Eq2.7}. Notice that the coefficients $\mathbf{F}$ and $\hat{\sigma}$ satisfy Assumption~\ref{AS1} and, moreover, if H\"{o}rmander's condition for parabolic hypoellipticity applies (see Remark~ref{R1} above). Then, the transition density exits and it satisfies the following Fokker-Planck equation
\begin{align}
\partial_{t} p(t, T, \mathbf{x},\mathbf{y}) + \mathcal{L}_{t,\mathbf{x}} p(t, T,\mathbf{x}, \mathbf{y}) &= 0, \quad 0 \le t < T, \quad \mathbf{x}, \mathbf{y} \in \mathbb{R}^3 \notag\\
p(T, T, \mathbf{x}, \mathbf{y}) &= \delta_{\mathbf{y}}(\mathbf{x}), \quad \mathbf{x}, \mathbf{y} \in \mathbb{R}^3, \label{Eq2.11}
\end{align}
where $\mathcal{L}_{t,\mathbf{x}} = \frac{1}{2} \operatorname{tr}\bigl \{a(t,\mathbf{x}) D_{x_1}^2 \bigr\} + \sum\nolimits_{i=1}^3 f_{i}(t,\mathbf{x}) D_{x_i}$ is the infinitesimal generator of $\mathbf{X}(t)$ at time $t$. On the other hand, starting from the above Fokker-Planck equation, if we replace $p(t, T, \mathbf{x},\mathbf{y})$ with $-\ln p(t, T, \mathbf{x}, \mathbf{y})$ (which is also known as the {\it Fleming's logarithmic transform}). Then, we should be able to establish a connection between the probabilistic representation for the transition density, i.e., $-\ln p(t, T, \mathbf{x}, \mathbf{y})$, with $0 \le t < T$, and that of the value function of a certain optimal stochastic control problem (e.g., see \cite{FleSh85} or \cite{She91} for a similar argument). 

In order to make our formulation mathematically more appealing, let us first fix the arrival point $(T, \mathbf{y}_0)$ and approximate the Dirac boundary condition in Equation~\eqref{Eq2.11} by a regular function, where we introduce a family of sequences $(\eta_{\varepsilon})_{\varepsilon > 0}$ of mollifier on $\mathbb{R}$ and weakly converges to the Dirac mass $\delta_{\mathbf{y}_0}$ at $\mathbf{y}_0$. Here, we also assume that $(\eta_{\varepsilon})_{\varepsilon > 0}$ are positive on $\mathbb{R}$. Then, we can approximate the transition density function by setting, for all $\varepsilon > 0$ and $(t, \mathbf{x}) \in [0, T - \varepsilon] \times \mathbb{R}^3$, $\upsilon_{\varepsilon}(t, \mathbf{x}) = \mathbb{E}_{t,\mathbf{x}} \bigl\{\eta_{\varepsilon}\bigl(\mathbf{X}(T - \varepsilon)\bigr) \bigr\}$.\footnote{Notice that the expectation $\mathbb{E}_{t,\mathbf{x}} \bigl\{\cdot\bigr\}$ is conditioned on $(t,\mathbf{x}) \in [0,  T- \varepsilon) \times \mathbb{R}^3$.} Note that, since the coefficients in Equation~\eqref{Eq2.11} are smooth, then such a solution is interpreted in the classical sense of Cauchy problem $\partial_t \upsilon_{\varepsilon}(t, \mathbf{x}) + \mathcal{L}_{t,\mathbf{x}} \upsilon_{\varepsilon}(t, \mathbf{x}) = 0$, for $0 \le t < T -\varepsilon$ and $\mathbf{x} \in \mathbb{R}^3$, with boundary condition $\upsilon_{\varepsilon}(t, \mathbf{x}) = \eta_{\varepsilon}\bigl(\mathbf{X}(T - \varepsilon)\bigr)$. Then, since $p$ is continuous away from the boundary, from the localization property of $(\eta_{\varepsilon})_{\varepsilon > 0}$, we can show that
\begin{align*}
\lim_{\varepsilon \rightarrow 0} \upsilon_{\varepsilon}(0, \mathbf{x}) &= \lim_{\varepsilon \rightarrow 0} \mathbb{E}_{0,\mathbf{x}} \bigl\{\eta_{\varepsilon}\bigl(\mathbf{X}(T - \varepsilon)\bigr) \bigr\}\\
                                                                                                 &= \lim_{\varepsilon \rightarrow 0} \int_{\mathbb{R}^3} \eta_{\varepsilon}\bigl(\mathbf{y}) p(0, T-\varepsilon,\mathbf{x}, \mathbf{y}) d\mathbf{y} \\
                                                                                                 &= p(0, T-\varepsilon,\mathbf{x}, \mathbf{y}_{0}).
\end{align*}
Moreover, for all $ t \in [0, T - \varepsilon]$ and $\mathbf{x} \in \mathbb{R}^3$, if we set $J_{\varepsilon}(t, \mathbf{x}) = -\ln \upsilon_{\varepsilon}(t, \mathbf{x})$, with $\upsilon_{\varepsilon}(t, \mathbf{x}) > 0$. Then, we can see that $J_{\varepsilon}(t, \mathbf{x})$ satisfies the following nonlinear parabolic partial differential equation (PDE)\footnote{$\langle \cdot, \cdot \rangle$ stands for the inner product.}
\begin{align}
\partial_{t} J_{\varepsilon}(t, \mathbf{x}) + \mathcal{L}_{t,\mathbf{x}} J_{\varepsilon}(t, \mathbf{x}) - \frac{1}{2} \bigl\langle a(s,\mathbf{x}) D_{x_1}J_{\varepsilon}(t, \mathbf{x}), D_{x_1} J_{\varepsilon}(t, \mathbf{x}) \bigr\rangle = 0, \label{Eq2.12}
\end{align}
with a boundary condition of $J_{\varepsilon}(T - \varepsilon, \mathbf{x}) = -\ln \eta_{\varepsilon}(\mathbf{x})$. Note that the above parabolic PDE can be rewritten as follows
\begin{align}
\partial_{t} J_{\varepsilon}(t, \mathbf{x}) &+ \mathcal{L}_{t,\mathbf{x}} J_{\varepsilon}(t, \mathbf{x}) \notag \\
& \quad + \inf_{u \in \mathbb{R}} \Bigl\{ \bigl\langle u,\, D_{x_1} J_{\varepsilon}(t, \mathbf{x}) \bigr\rangle + \frac{1}{2} \bigl\langle a^{-1}(t,\mathbf{x}) u,\, u \bigr\rangle \Bigr\} = 0. \label{Eq2.13}
\end{align}
Then, for $(t, \mathbf{x}) \in [0, T - \varepsilon) \times \mathbb{R}^3$, the infimum in the above equation, i.e., Equation~\eqref{Eq2.13}, is achieved, when $u(t, \mathbf{x})=-a(t,\mathbf{x}) D_{x_1}J_{\varepsilon}(t, \mathbf{x})$, which exactly gives us the relation in Equation~\eqref{Eq2.12}.

Next, let us denote by $\mathcal{U}_{T - \varepsilon}$ the set of $\mathbb{R}$-valued progressively measurable processes $\bigl(u(t)\bigr)_{0 \le t \le T - \varepsilon}$ (i.e., a family of nonanticipative processes, for all $t > s$, $(W(t)-W(s))$ is independent of $u(r)$ for $r \le s$) such that
\begin{align*}
\mathbb{E} \int_{0}^{T - \varepsilon} \vert u(t)\vert^2 dt < \infty.
\end{align*}
Moreover, from Equation~\eqref{Eq2.13}, for $t \in [0, T - \varepsilon]$, with $u(t) \in \mathcal{U}_{T - \varepsilon}$, we can write $J_{\varepsilon}(t, \mathbf{x})$ as the value function for the following optimal stochastic control problem
\begin{align}
J_{\varepsilon}(0, \mathbf{x}) = \inf_{u(t) \in \mathcal{U}_{T - \varepsilon}} \mathbb{E}_{0, \mathbf{x}} \Bigl\{\frac{1}{2} \int_{0}^{T - \varepsilon} \bigl\langle a^{-1}(t,\mathbf{X}_{0,\mathbf{x}}^{u}(t)) u(t),\, u(t)\bigr\rangle dt \notag \\
- \ln \eta_{\varepsilon}(\mathbf{X}_{0,\mathbf{x}}^{u}(T - \varepsilon)) \Bigr\}, \label{Eq2.14}
\end{align}
which is associated with a controlled version of SDE in Equation~\eqref{Eq2.7} (with the corresponding controlled-diffusion process $\bigl(\mathbf{X}_{0,\mathbf{x}}^{u}(t) \bigr)_{t \ge 0}$), i.e., 
\begin{align}
d \mathbf{X}_{0,\mathbf{x}}^{u}(t) = \bigl[\mathbf{F} (t, \mathbf{X}_{0,\mathbf{x}}^{u}(t)) + B u(t) \bigr]dt + B \hat{\sigma}(t, \mathbf{X}_{0,\mathbf{x}}^{u}(t)) dW(t), \,\, \mathbf{X}_{0,\mathbf{x}}^{u}(0) = \mathbf{x}. \label{Eq2.15}
\end{align}
Note that if we expand Equation~\eqref{Eq2.14} using It\^{o}'s formula, then, from Equation~\eqref{Eq2.12}, we obtain the following
\begin{align}
& dJ_{\varepsilon}(t, \mathbf{X}_{0,\mathbf{x}}^{u}(t)) = \Bigl \{\partial_{t} J_{\varepsilon}(t,\mathbf{X}_{0,\mathbf{x}}^{u}(t)) + \mathcal{L}_{t,\mathbf{x}} J_{\varepsilon}(t, \mathbf{x}) +  \bigl\langle D_{x_1}J_{\varepsilon}(t,\mathbf{X}_{0,\mathbf{x}}^{u}(t)), u(t) \bigr\rangle \Bigr \} dt \notag \\
                            &\quad\quad\quad \quad  + \bigl\langle D_{x_1} J_{\varepsilon}(t,\mathbf{X}_{0,\mathbf{x}}^{u}(t)), \hat{\sigma}(t, \mathbf{X}_{0,\mathbf{x}}^{u}(t)) dW(t) \bigr\rangle \notag \\
                             & \quad \quad = \left \{\frac{1}{2} \bigl\langle D_{x_1}J_{\varepsilon}(t,\mathbf{X}_{0,\mathbf{x}}^{u}(t))u^{\ast}(t), u^{\ast}(t) \bigr\rangle - \bigl\langle D_{x_1}J_{\varepsilon}(t,\mathbf{X}_{0,\mathbf{x}}^{u}(t))u(t), u^{\ast}(t) \bigr\rangle   \right \} dt \notag \\
                            &\quad \quad\quad \quad  + \bigl\langle \hat{\sigma}^{-1}(t,\mathbf{X}_{0,\mathbf{x}}^{u}(t)) u^{\ast}(t), dW(t) \bigr\rangle \notag \\
                             & \quad \quad = \frac{1}{2} \left \{ \bigl\vert \hat{\sigma}^{-1}(t,\mathbf{X}_{0,\mathbf{x}}^{u}(t)) \bigl(u^{\ast}(t) - u(t)\bigr) \bigr\vert^2 - \bigl \vert \hat{\sigma}^{-1}(t,\mathbf{X}_{0,\mathbf{x}}^{u}(t)) u(t) \bigr\vert^2   \right \} dt \notag \\
                            &\quad \quad\quad \quad + \bigl\langle \hat{\sigma}^{-1}(t,\mathbf{X}_{0,\mathbf{x}}^{u}(t)) u^{\ast}(t), dW(t) \bigr\rangle, \label{Eq2.16}
\end{align}
where $u^{\ast} = - a(t,\mathbf{X}_{0,\mathbf{x}}^{u}(t)) D_{x_1}J_{\varepsilon}(t,\mathbf{X}_{0,\mathbf{x}}^{u}(t))$, so that

\begin{align}
J_{\varepsilon}(0, \mathbf{x}) =& -\ln \eta_{\varepsilon}(\mathbf{X}_{0,\mathbf{x}}^{u}(t)) + \frac{1}{2} \int_{0}^{T - \varepsilon} \bigl\langle D_{x_1}J_{\varepsilon}(t,\mathbf{X}_{0,\mathbf{x}}^{u}(t))u(t), u(t) \bigr\rangle dt \notag \\
                                        & \quad - \frac{1}{2} \int_{0}^{T - \varepsilon}  \bigl\vert \hat{\sigma}^{-1}(t,\mathbf{X}_{0,\mathbf{x}}^{u}(t)) \bigl(u^{\ast}(t) - u(t)\bigr) \bigr\vert^2 dt \notag \\
                                        & \quad + \int_{0}^{T - \varepsilon} \bigl\langle \hat{\sigma}^{-1}(t,\mathbf{X}_{0,\mathbf{x}}^{u}(t)) u^{\ast}(t), dW(t) \bigr\rangle. \label{Eq2.18}
\end{align}

Before concluding this subsection, let us briefly discuss the deterministic version of the above optimal stochastic control problem in Equations~\eqref{Eq2.14} and \eqref{Eq2.15}. To this end, let $\bigl(\pmb{\phi}(t)\bigr)_{0 \le t \le T}$ be the solution of the following deterministic control problem 
\begin{align}
\dot{\pmb{\phi}}(t) = \mathbf{F} (t, \pmb{\phi}(t)) + B \varphi(t), \quad \pmb{\phi}(0) = \mathbf{x}_0, \label{Eq2.19}
\end{align}
where $\bigl(\varphi(t)\bigr)_{0 \le t \le T}$ is a deterministic control from the space $L^2([0,T], \mathbb{R})$ and it is chosen so as to force $\pmb{\phi}$ to hit $\mathbf{y}_0$ at time $T$, i.e., $\pmb{\phi}(T) = \mathbf{y}_0$, with a minimum quadratic cost that is exactly as in Equation~\eqref{Eq2.14}, where $\bigl(u(t)\bigr)_{0 \le t \le T - \varepsilon}$ is chosen to make $\mathbf{X}_{0,\mathbf{x}_0}^{u}(T - \varepsilon)$ tends to $\mathbf{y}_0$, when $\varepsilon \rightarrow 0$. In other words, we choose $\bigl(\varphi(t)\bigr)_{0 \le t \le T}$ from $L^2([0,T], \mathbb{R})$ as an optimal control for the following minimization problem
\begin{align}
I(0, T, \mathbf{x}_0, \mathbf{y}_0) = \inf \Bigl\{\frac{1}{2} \int_{0}^{T} \bigl \vert \varphi(t) \bigr\vert^2 dt \,\Bigl \vert \,\, \pmb{\phi}(0) = \mathbf{x}_0, \,\, \pmb{\phi}(T) = \mathbf{y}_0 \Bigr\}. \label{Eq2.20}
\end{align}
Moreover, the functional $I(0, T, \mathbf{x}_0, \mathbf{y}_0)$ is finite, since the dynamical system in Equation~\eqref{Eq2.19} is assumed to be controllable (cf. \cite{Ell73}), and such a functional also satisfies the same order as that of $\bigl\vert T^{1/2}\Gamma^{-1}(T) \bigl(\Theta(T, \mathbf{x}_0) - \mathbf{y}_0 \bigr)\bigr\vert^2$ (see also Section~\ref{S3} for an intrinsic time-scaling property), where $\Theta(T, \mathbf{x}_0)$ denotes the deterministic solution to the dynamical system in Equation~\eqref{Eq2.2} at time $T$ (cf. Equation~\eqref{Eq2.2c} below).

\begin{remark} \label{R2}
Note that $I(0, T, \mathbf{x}, \mathbf{y})$ is known as the action functional in large deviations theory that provides, in short-time, a natural link between the deterministic control problem in Equations~\eqref{Eq2.19}-\eqref{Eq2.20} and that of the transition density function of Equation~\eqref{Eq2.11} (e.g., see \cite{FreWe84} or \cite{VenFre70} for additional discussions).
\end{remark}

In the following section, i.e., in Section~\ref{S3}, using the solution of the ``{\it functional equation}" for the value function in Equation~\eqref{Eq2.14} together with the limit for $J_{\varepsilon}(0, \mathbf{x})$, when $\varepsilon \rightarrow 0$ (i.e., $\lim_{\varepsilon \rightarrow 0} J_{\varepsilon}(0, \mathbf{x}) = -\ln p(0, T, \mathbf{x}, \mathbf{y})$), we provide bounds on the transition density function $p(0, T, \mathbf{x}, \mathbf{y})$.

\section{Main result} \label{S3}
In this section, we provide bounds on the solution of the transition density function for the Fokker-Planck equation in Equation~\eqref{Eq2.11} that corresponds to the opioid epidemic dynamical model, when the random perturbation enters only through the dynamics of the susceptible group in the compartmental model (i.e., the SDE in Equation~\eqref{Eq2.6} or \eqref{Eq2.7}). Here, the proof for such bounds basically relies on the interpretation of the solution for the density function as the value function of some optimal stochastic control problem (see also Equations~\eqref{Eq2.14} and \eqref{Eq2.15}).

For any $t \in [0, T]$, let $\Theta(t, \mathbf{x}) = (\theta_1(t), \theta_2(t), \theta_3(t))$ be the solution of the system of (deterministic) ordinary differential equations in Equation~\eqref{Eq2.2} (see also Equation~\eqref{Eq2.2b}), i.e., 
\begin{eqnarray}
\left.\begin{array}{l}
\dot{\theta}_1(t) = f_1(t, \theta_1, \theta_2, \theta_3)\\
\dot{\theta}_2(t) = f_2(t, \theta_1, \theta_2, \theta_3) \\
\dot{\theta}_3(t) = f_3(t, \theta_2, \theta_3)
\end{array}\right\}  \,\, \rightsquigarrow \,\, \Theta(t, \mathbf{x}) = \bigl(\theta_1(t), \theta_2(t), \theta_3(t) \bigr), \label{Eq2.2c}
\end{eqnarray}
with an initial condition of $\Theta(0, \mathbf{x}) = (\theta_1(0), \theta_2(0), \theta_3(0)) = \mathbf{x} \in \mathbb{R}^3$, starting from time zero. Then, we have the following result useful for characterizing the solution of the transition density function in Equation~\eqref{Eq2.11}.
\begin{proposition} \label{P3.1}
Suppose that the statements in Assumption~\ref{AS1} hold, then, for any time $t \in [0, T]$, the law $\mathbf{X}(t)$, i.e., the solution of the SDE in Equation~\eqref{Eq2.7} starting from an initial condition $\mathbf{x} \in \mathbb{R}^3$ and an initial time zero, admits a transition density function $\mathbf{y} \in \mathbb{R}^3 \mapsto p(0, t, \mathbf{x}, \mathbf{y})$ that satisfies the Fokker-Planck equation in Equation~\eqref{Eq2.11}. Moreover, for any $T > 0$, there exists a constant $\mathbf{C} > 1$, such that, for any $t \in [0, T]$,
\begin{align}
\mathbf{C}^{-1} t^{-9/2} \exp \Bigl( -\mathbf{C} t \Gamma(t) \bigl\vert \Theta(t, \mathbf{x}) - & \mathbf{y} \bigr \vert^2 \Bigr) \notag \\
              \le p(0, t, \mathbf{x}, \mathbf{y}) &  \le  \mathbf{C} t^{-9/2} \exp \Bigl( -\mathbf{C} t \Gamma(t) \bigl\vert \Theta(t, \mathbf{x}) -  \mathbf{y} \bigr \vert^2 \Bigr),
\end{align}
where the intrinsic time-scaling matrix $\Gamma(t)$ is a $3 \times 3$ diagonal matrix which is given by $\Gamma(t) = T^{-1}\diag\bigl\{t,\, t^2,\, t^3\bigr\}$. \label{Eq3.1}
\end{proposition}

Before attempting to prove the above proposition, let us briefly discuss the time-scaling property with which the system dynamics in Equation~\eqref{Eq2.6} is evolving. Note that the random noise enters only through the dynamics of the susceptible group in Equation~\eqref{Eq2.2} and is then subsequently propagated to the other groups in the compartmental model. Moreover, from the weak H\"{o}rmander's condition, we also observe that the drift terms are responsible for propagating the noise through the system dynamics. As a result, the underlying constant $C$ in the above proposition in general depends on $T$, whereas the statement around Equation~\eqref{Eq3.1} is valid for any time $t \in [0, T]$. Hence, one possible way of addressing this issue is to use the intrinsic time-scaling property in Equation~\eqref{Eq2.6} and to pass from $t$ to $T$ (or passing from $T$ to $1$). For example, for a given $T > 0$, let us define the rescaled version of $\bigl(\widetilde{\mathbf{X}}(t)\bigr)_{t \ge 0}$ by setting
\begin{align*}
\widetilde{\mathbf{X}}(t) &= T^{1/2} \Gamma^{-1}(T) \mathbf{X}(t)\\
                               &= \bigl(T^{-1/2} X_1(Tt), T^{-3/2} X_2(Tt), T^{-5/2}X_3(Tt)\bigr), \quad t \ge 0,
\end{align*}
where $\Gamma(T)$ represents the scaling matrix at time $T$. 

Then, the rescaled diffusion process $\bigl(\widetilde{\mathbf{X}}(t)\bigr)_{t \ge 0}$ satisfies the following
\begin{align}
d \widetilde{\mathbf{X}}(t) =& T^{-3/2} \Gamma^{-1}(T) \mathbf{F} (Tt, T^{-1/2} \Gamma(T) \widetilde{\mathbf{X}}(t)) dt \notag \\
                                      & \quad \quad \quad \quad + B \hat{\sigma}(Tt, T^{-1/2} \Gamma(T) \widetilde{\mathbf{X}}(t)) d\tilde{W}(t), \label{Eq3.2}
\end{align}
where $\widetilde{W}(t)$ stands for the rescaled Brownian motion, i.e., 
\begin{align*}
\bigl(\widetilde{W}(t)\bigr)_{t \ge 0} = \bigl(T^{1/2} W(Tt)\bigr)_{t \ge 0}.
\end{align*}
If we further denote by $\mathbf{y} \in \mathbb{R}^3 \mapsto \widetilde{p}(0, 1, \mathbf{x}, \mathbf{y})$ the corresponding transition density function at time $1$. Then, with $\det( T^{-1/2} \Gamma(T))=T^{9/2}$, we obtain the following
\begin{align}
\widetilde{p}(0, 1, \mathbf{x}, \mathbf{y}) = T^{9/2} p(0, T, T^{-1/2} \Gamma(T) \mathbf{x}, T^{-1/2} \Gamma(T)\mathbf{y}), \quad \forall \mathbf{x}, \mathbf{y} \in \mathbb{R}^3. \label{Eq3.3}
\end{align}
Equivalently, we can also obtain an estimate for the transition density function $p(0, T, \mathbf{x}, \mathbf{y})$ at time $T$ from an estimate at time $1$, i.e.,
\begin{align}
p(0, T, \mathbf{x}, \mathbf{y}) = T^{-9/2} \widetilde{p}(0, 1, T^{1/2} \Gamma^{-1}(T)\mathbf{x}, T^{1/2} \Gamma^{-1}(T) \mathbf{y}), \quad \forall \mathbf{x}, \mathbf{y} \in \mathbb{R}^3. \label{Eq3.4}
\end{align}
Moreover, the transition density function $\widetilde{p}(0, 1, \mathbf{x}, \mathbf{y})$, for all $\mathbf{x}, \mathbf{y} \in \mathbb{R}^3$, admits two-side bounds, with $\mathbf{C} \exp \bigl( -\mathbf{C} \bigl\vert \widetilde{\Theta}(t, \mathbf{x}) - \mathbf{y} \bigr \vert^2 \bigr)$, where $\bigl(\widetilde{\Theta}(t, \mathbf{x})\bigr)_{t \ge 0}$ represents the rescaled flow and is given by
\begin{align*}
\widetilde{\Theta}(t, \mathbf{x}) = T^{1/2} \Gamma^{-1}(T) \Theta(Tt, T^{-1/2} \Gamma(T)\mathbf{x}),
\end{align*}
for $(t, \mathbf{x}) \in [0, 1] \times \mathbb{R}^3$. Note that, in the proof part below, we use this intrinsic property to support our claim, i.e., for providing bounds on the solution of the transition density function to the Fokker-Planck equation in Equation~\eqref{Eq2.11} (see also \cite{Aro67} for additional discussion on bounds for fundamental solutions of parabolic equations).

{\em Sketch of the proof}.
The proof for the above proposition (i.e., Proposition~\ref{P3.1} which is an adaptation of \cite{DelM10}) involves finding bounds for $J_{\varepsilon}(0, \mathbf{x})$, uniformly in $\varepsilon > 0$. Note that $J_{\varepsilon}(0, \mathbf{x})$ is the value function for the minimization problem in Equation~\eqref{Eq2.18}. Hence, the whole problem then is reduced in finding an admissible control process $\bigl(u(t)\bigr)_{0 \le t \le T - \varepsilon}$ from the set $\mathcal{U}_{T - \varepsilon}$, for each $\varepsilon > 0$. 

In the following, we proceed as follows: we first perform a standard linearization of the controlled SDE in Equation~\eqref{Eq2.15} around the deterministic curve solution $\bigl(\pmb{\phi}(t)\bigr)_{0 \le t \le T}$ of Equation~\eqref{Eq2.19}, i.e., expanding $\mathbf{F} (t, \mathbf{X}_{0,\mathbf{x}}^{u}(t))$ as follows
\begin{align*}
\mathbf{F} (t, \mathbf{X}_{0,\mathbf{x}}^{u}(t)) = \mathbf{F} (t, \pmb{\phi}(t)) &+ D_{\mathbf{x}}\mathbf{F} (t, \pmb{\phi}(t)) \bigl(\mathbf{X}_{0,\mathbf{x}}^{u}(t) - \pmb{\phi}(t) \bigr) \\
& \quad\quad + o\bigl(\bigl\vert \mathbf{X}_{0,\mathbf{x}}^{u}(t) - \pmb{\phi}(t) \bigr\vert \bigr),
\end{align*}
where $D_{\mathbf{x}}$ denotes the space derivative, i.e., w.r.t. $\mathbf{x}$, of $\mathbf{F}$ and we similarly approximate $\hat{\sigma}(t,\mathbf{X}_{0,\mathbf{x}}^{u}(t))$ by $\hat{\sigma}(t, \pmb{\phi}(t))$. Note that the curve $\bigl(\pmb{\phi}(t)\bigr)_{0 \le t \le T}$ is deterministic and the following mapping  
\begin{align*}
 \mathbb{R}^3 \ni \mathbf{z} \mapsto \mathbf{F} (t, \pmb{\phi}(t)) + D_{\mathbf{x}}\mathbf{F} (t, \pmb{\phi}(t)) \bigl(\mathbf{z} - \pmb{\phi}(t) \bigr)
\end{align*}
is affine. Noting the approximate diffusion coefficient $\hat{\sigma}(t, \pmb{\phi}(t))$, for $0 \le t \le T$, is deterministic. Then, from Equations~\eqref{Eq2.15} and \eqref{Eq2.19}, we also observe that $\bigl(\mathbf{X}_{0,\mathbf{x}}^{u}(t) - \pmb{\phi}(t)\bigr)$, for $0 \le t \le T - \varepsilon$ satisfies the following linearized system\footnote{Higher order terms are assumed to be well controlled.}
\begin{align}
d\bigl(\mathbf{X}_{0,\mathbf{x}}^{u}(t) - \pmb{\phi}(t)\bigr) \approx \bigl(D_{\mathbf{x}}\mathbf{F} (t, \pmb{\phi}(t)) \bigl(\mathbf{X}_{0,\mathbf{x}}^{u}(t) - \pmb{\phi}(t)\bigr) &+ B \bigl(u(t) - \varphi(t) \bigr) \bigr)dt \\
& + B \hat{\sigma}(t, \pmb{\phi}(t)) dW(t), \label{Eq3.5}
\end{align}
where the initial point $\bigl(\mathbf{X}_{0,\mathbf{x}}^{u}(0) - \pmb{\phi}(0)\bigr)$ is zero and that of the final point $\bigl(\mathbf{X}_{0,\mathbf{x}}^{u}(T - \varepsilon) - \pmb{\phi}(T - \varepsilon)\bigr)$ is also expected to be close to zero. In some sense, we reduced to the linear problem, with zero boundary conditions and $\bigl(u(t) - \varphi(t) \bigr)$, for $0 \le t \le T - \varepsilon$ as an admissible control (or equivalently determining an appropriate control $u(t) \in \mathcal{U}_{T - \varepsilon}$).

Note that, for the above linearized system (i.e., Equation~\eqref{Eq3.5}) with constant diffusion term, the corresponding optimal stochastic control problem (cf. Equations~\eqref{Eq2.14} and \eqref{Eq2.15}) can be written in an explicit form. Moreover, $\bigl(\mathbf{X}_{0,\mathbf{x}}^{u}(t)\bigr)$ is a Gaussian process and, hence, the transition density function exists and has an explicit form that satisfies two-side exponential bounds of Equation~\eqref{Eq3.1} (see also \cite[Theorem~1]{Aro67}). This completes the proof of Proposition~\ref{P3.1}. \qed

\section{On the attainable exit probability} \label{S4} 
In this section, for a fixed $T > 0$, we provide an estimate for the attainable exit probability with which the diffusion process $\mathbf{X}(t)$ exits from a given bounded open domain $D \subset \mathbb{R}^3$ with smooth boundary $\partial D$ (i.e., $\partial D$ is a manifold of class $C^2$), that is,
\begin{align*}
q(\mathbf{x}) = \mathbb{P}_{0,\mathbf{x}}\bigl\{\tau \le T \bigr\},
\end{align*}
where $\tau$ is the first exit-time for the diffusion process $\mathbf{X}(t)$ from the domain $D$, i.e., $\tau = \inf \bigl\{ t > 0 \, \bigl\vert \, \mathbf{X}(t) \notin D \bigr\}$. Note that the exit point distribution is also intimately connected with a family of Dirichlet problems on $\Omega$ that is associated with a nondegenerate diffusion process $\mathbf{X}^{\epsilon}(t)$ (see Equation~\eqref{Eq4.2} below) as the limiting case, when $\epsilon \rightarrow 0$ (e.g., see \cite{BefA15a} for related discussions on the attainable exit probabilities for the diffusion processes).

\begin{proposition} \label{P4.1}
Let $\psi (\mathbf{x})$ be a continuous function on $\partial D$. Then, the attainable exit probability $q(\mathbf{x})$ with which the diffusion process $\mathbf{X}(t)$ exits from the domain $D$ is a smooth solution to the following Dirichlet problem
\begin{align}
\left.\begin{array}{c}
\mathcal{L} \upsilon(\mathbf{x}) = 0 \quad \text{in} \quad D  \\
 \upsilon(\mathbf{x}) = \mathbb{E}_{0,\mathbf{x}} \bigl\{ \psi \bigl(\mathbf{X}(\tau) \bigr)\bigr\} \quad \text{on} \quad  \partial D
\end{array}\right\}.  \label{Eq4.1}
\end{align}
Moreover, it is a continuous function on $D$.
\end{proposition}

In order to prove the above proposition (which an adaptation of the proof of \cite[Proposition~3.5]{BefA15a}), we will consider a nondegenerate diffusion process $\mathbf{X}^{\epsilon}(t) = (X_1^{\epsilon}(t), X_2^{\epsilon}(t), X_3^{\epsilon}(t))$, i.e.,
\begin{eqnarray}
\left.\begin{array}{l}
dX_1^{\epsilon}(t) = f_1(t, X_1^{\epsilon}(t), X_2^{\epsilon}(t), X_3^{\epsilon}(t)) dt + \hat{\sigma}(t, X_1^{\epsilon}(t), X_2(t), X_3^{\epsilon}(t)) dW(t)\\
dX_2^{\epsilon}(t) = f_2(t, X_1^{\epsilon}(t), X_2^{\epsilon}(t), X_3^{\epsilon}(t)) dt + \sqrt{\epsilon} dV(t)\\
dX_3^{\epsilon}(t) = f_3(t, X_2^{\epsilon}(t), X_3^{\epsilon}(t)) dt + \sqrt{\epsilon} dV(t)
\end{array}\right\} \label{Eq4.2}
\end{eqnarray}
where $V(t)$ (with $V(0)=0$) is a $1$-dimensional brownian motion independent to $W(t)$. Then, we relate the exit probability of this diffusion process $\mathbf{X}^{\epsilon}(t)$ with that of the Dirichlet problem in Equation~\eqref{Eq4.1} in the limiting case, when $\epsilon \rightarrow 0$.

Let us define the following notations that will be useful in the sequel.
\begin{align*}
\begin{array}{c}
    \theta = \tau \wedge T, \quad\quad  \theta^{\epsilon} = \tau^{\epsilon} \wedge T,\\ \\
   \bigl\Vert \mathbf{X}^{\epsilon} - \mathbf{X} \bigr\Vert_t = \sup\limits_{0 \le r \le t} \bigl\vert \mathbf{X}^{\epsilon}(r) - \mathbf{X}(r) \bigr\vert,
\end{array}
\end{align*}
where $\tau^{\epsilon}$ is the first exit-time for the diffusion process $\mathbf{X}^{\epsilon}(t)$ from the domain $D$, i.e., $\tau^{\epsilon} = \inf \bigl\{ t > 0 \, \bigl\vert \, \mathbf{X}^{\epsilon}(t) \notin D \bigr\}$. Moreover, let $n(\mathbf{x})$ be an outer normal vector to the boundary of $D$ and we further denote by $\partial D^{+}$ the set of points $\mathbf{x} \in \mathbb{R}^3$, with $\mathbf{x} \in \partial D$, such that $\langle f_1(t,\mathbf{x}), n(\mathbf{x})\rangle$ is positive.

\begin{proof}
Note that, from Assumption~\ref{AS1}(b), it is sufficient to show that $q(\mathbf{x})$ is a smooth solution (almost everywhere in $D$ with respect to Lebesgue measure) to the Dirichlet problem in Equation~\eqref{Eq4.1}.

Consider the following infinitesimal generator that corresponds to the above nondegenerate diffusion process $\mathbf{X}^{\epsilon}(t)$, with $\epsilon \ll 1$,
\begin{align}
 \mathcal{L} \upsilon(\mathbf{x})  + \frac{\epsilon}{2} \sum\nolimits_{i=2}^3 \triangle_{x_{i}} \upsilon(\mathbf{x})  = 0 \quad \text{in} \quad D,  \label{Eq4.3}
\end{align}
where $\triangle_{x_{i}}$ is the Laplace operator in the variable $x_{i}$ and $\mathcal{L}$ is the infinitesimal generator in Equation~\eqref{Eq2.11}.

Note that the infinitesimal generator $\mathcal{L}$ is a uniformly parabolic equation and, from Assumption~\ref{AS1}(b), its solution satisfies the following boundary condition
\begin{align}
 \upsilon(\mathbf{x}) = \psi (\mathbf{x}) \quad \text{on} \quad \partial D,  \label{Eq4.4}
\end{align}
where\footnote{Here, $\mathbb{E}_{0,\mathbf{x}}^{\epsilon} \bigl\{\cdot\bigr\}$ is associated with the diffusion process $\mathbf{X}^{\epsilon}(t)$.}
\begin{align}
 \upsilon (\mathbf{x} ) = \mathbb{E}_{0,\mathbf{x}}^{\epsilon} \bigl\{ \psi (\mathbf{X}^{\epsilon}(\theta^{\epsilon}) ) \bigr\},  \label{Eq4.5}
\end{align}
with $\theta^{\epsilon} = \tau^{\epsilon} \wedge T$.

In particular, let $\psi_k$, with $k=1, 2,\ldots$, be a sequence of bounded functions that are continuous on $\partial D$ and satisfying the following conditions
\begin{align*}
 \psi_k \bigl(\mathbf{X}^{\epsilon}(t)\bigr) = \left\{\begin{array}{l l}
1  \quad &\text{if} \quad \mathbf{X}^{\epsilon}(t) \in \partial D^{+}\\
0 \quad &\text{if}  \quad \mathbf{X}^{\epsilon}(t) \in D \quad \text{and} \quad \varrho\bigl(\mathbf{X}^{\epsilon}(t), \partial D\bigr) > \frac{1}{k}
\end{array}\right. 
\end{align*}
and
\begin{align*}
0 \le \psi_k \bigl(\mathbf{X}^{\epsilon}(t)\bigr) \le 1 \quad & \text{if} \quad  \mathbf{X}^{\epsilon}(t) \in D \quad \text{and} \quad  \varrho\bigl(\mathbf{X}^{\epsilon}(t), \partial D\bigr) \le \tfrac{1}{k}.
\end{align*}
Moreover, the bounded functions further satisfy the following
\begin{align}
\bigl\vert \psi_k - \psi_m \bigr\vert \rightarrow 0 \quad \text{as} \quad k, m \rightarrow \infty  \label{Eq4.6}
\end{align}
uniformly on any compact subset of $\bar{D}$. Then, with $\psi = \psi_k$,
\begin{align*}
 \upsilon_k (\mathbf{x}) = \mathbb{E}_{0,\mathbf{x}} \bigl\{ \psi_k \bigl(\mathbf{X}^{\epsilon}(\theta^{\epsilon}) \bigr) \bigr\}
\end{align*}
satisfies Equations~\eqref{Eq4.5} and \eqref{Eq4.6}. Then, from the continuity of $\psi_k$ (cf. \cite[Lemma~2.5]{BefA15a}, Parts~(i)-(iii)) and the Lebesgue's dominated convergence theorem (e.g., see \cite[Chapter~4]{Roy88}), we have the following
\begin{align}
 \upsilon_k (\mathbf{x}) \rightarrow \underbrace{\mathbb{E}_{0,\mathbf{x}}^{\epsilon} \bigl\{\psi_k \bigl(\mathbf{X}(\theta^{\epsilon}\bigr) \bigr)\bigr\}}_{\triangleq  q_k(\mathbf{x})}, \label{Eq4.7}
\end{align}
where $\theta^{\epsilon} \rightarrow \theta$ and $\bigl\Vert \mathbf{X}^{\epsilon} - \mathbf{X} \bigr\Vert_T = \sup\limits_{0 \le r \le T} \bigl\vert \mathbf{X}^{\epsilon}(r) - \mathbf{X}(r) \bigr\vert \rightarrow 0$, as $\epsilon \rightarrow 0$; and $\mathbf{X}(t)$ is a solution to Equation~\eqref{Eq2.7}, when $\epsilon=0$, with an initial condition $\mathbf{X}(0) = \mathbf{x}$. 

Note that $q_k (\mathbf{x})$ satisfies Equation~\eqref{Eq4.3}, with $\upsilon (\mathbf{x}) = \upsilon_k (\mathbf{x})$, and, in addition, it is a distribution solution to the Dirichlet problem in Equation~\eqref{Eq4.1}, i.e., with any test function $\phi \in C_0^{\infty}(D)$,
\begin{align*}
\int_{D} \mathcal{L} \phi (\mathbf{x}) q_k (\mathbf{x}) d D &= \lim_{\epsilon \rightarrow 0} \int_{D} \Bigl(\mathcal{L} \phi(\mathbf{x}) + \frac{\epsilon}{2} \sum\nolimits_{i=2}^3  \triangle_{x_{i}} \phi (\mathbf{x})  \Bigr) \upsilon_k(\mathbf{x}) d D, \\
  &=0,
\end{align*}
where $dD$ denotes the Lebesgue measure on $\mathbb{R}^3$.

Finally, notice that (cf. Equation~\eqref{Eq4.7})
\begin{align*}
q (\mathbf{x}) = \lim_{k \rightarrow \infty} q_k (\mathbf{x}),
\end{align*}
almost everywhere in $D$. Moreover, from Assumption~\ref{AS1}(b) (i.e., the hypoellipticity), $q(\mathbf{x})$ is a smooth solution to Equation~\eqref{Eq4.1} (almost everywhere) in $D$ and continuous on the boundary of $D$. This completes the proof of Proposition~\ref{P4.1}. \qed
\end{proof}

\begin{remark} \label{R3}
Here, it is worth remarking that Propositions~\ref{P3.1} and \ref{P4.1} are useful for selecting the most appropriate admissible strategies that confines the diffusion process $\mathbf{X}(t)$ to the prescribed domain $D$ for a longer duration. Moreover, such an admissible strategy also depends on the state information (i.e.,  $u = - a(t,\mathbf{x}) D_{x_1}J_{\varepsilon}(t,\mathbf{x})$) and it only enters through the susceptible group in the compartmental model (cf. Equation~\eqref{Eq2.15}).
\end{remark}

\section{Concluding remarks} \label{S5}
In this paper, we have considered two seemingly related problems pertaining to an opioid epidemic dynamical model with random perturbation. In the first problem, we have provided two-sided bounds on the solution of the transition density function for the Fokker-Planck equation that is associated with the underlying opioid epidemic dynamical model, when the random perturbation enters only through the dynamics of the susceptible group in the compartmental model. In particular, we have argued that such bounds can be obtained based on a precise interpretation of the transition density function as a value function for a certain stochastic control problem. In the second problem, we have also provided an estimate on the attainable exit probability with which the solution for the randomly perturbed opioid epidemic dynamical model exits from a given bounded open domain during a certain time interval. Note that such qualitative information on the first exit-time as well as two-sided bounds on the transition density function are useful for developing effective and fact-informed intervention strategies that primarily aim at curbing opioid epidemics or assisting in interpreting outcome results from opioid-related policies.


\begin{thebibliography}{99}

\bibitem{r1}
Centers for Disease Control and Prevention. {\em Understanding the epidemic}.\\
Available at: https://www.cdc.gov/drugoverdose/epidemic/index.html. Accessed April 12, 2018. 

\bibitem{r2}
Volkow N.D., McLellan A.T.
{\em Opioid abuse in chronic pain -- misconceptions and mitigation strategies}.
N. Engl. J Med. 374, 1253--1263 (2016).

\bibitem{r3}
Frieden, T.R., Houry D.
{\em Reducing the risks of relief -- the CDC opioid-prescribing guideline}.
N. Engl. J. Med. 374, 1501--1504 (2016).

\bibitem{r4}
Meldrum, M.L.
{\em The ongoing opioid prescription epidemic: historical context}. 
Am. J. Public Health. 106, 1365--1366 (2016).

\bibitem{r6}
Battista, N.A., Pearcy, L.B., Strickland, W.C.
{\em Modeling the opioid epidemic}.
Preprint \href{https://arxiv.org/abs/1711.03658}{arXiv:1711.03658} \href{http://arxiv.org/archive/q-bio.PE}{[q-bio.PE]}, 27 pages, November 2017

\bibitem{r7}
Njagarah, H.J.B., Nyabadza, F.
{\em Modeling the impact of rehabilitation, amelioration and relapse on the prevalence of drug epidemics}. 
Journal of Biological Systems, 21(1), Article ID 1350001 (2013).

\bibitem{r8}
Kermack, W.O., McKendrick, A.G.
{\em A contribution to the mathematical theory of epidemics}.
Proc. Roy. Soc. Lond. A. 115, 700--721 (1927).

\bibitem{r9}
Samanta, G.
{\em Dynamic behaviour for a non autonomous heroin epidemic with time delay}.
J. Appl. Math. Comput. 35, 161--178 (2011).

\bibitem{r10}
Huang, G., Liu, A.
{\em A note on global stability for a heroin epidemic model with distributed delay}.
Appl. Math. Letters. 26, 687--691 (2013).

\bibitem{r11}
Bin, F., Xuezhi, L., Maia, M., Liming, C.
{\em Global stability for a heroin model with age-dependent susceptibility}.
J. Sys. Sci. and Complex. 28, 1243--1257 (2015).

\bibitem{r12}
Ma, M., Liu, S., Li, J.
{\em Bifurcation of a heroin model with nonlinear incidence rate}.
Nonlinear Dynamics. 88, 555--565 (2017).

\bibitem{r13}
White, E., Comiskey, C.
{\em Heroin epidemics, treatment and ode modelling}.
Math. Biosci. 208, 312--324 (2007).

\bibitem{CasJR97}
Castillo, G.C., Jordan, S.G., Rodriguez, A.H. 
{\em Mathematical models for the dynamics of tobacco use, recovery and relapse}.
Vol. {\bf 3}, Technical Report Series, BU-1505-M, Department of Biometrics, Cornell University (1997).

\bibitem{r5}
Ryan Haight Online Consumer Protection Act of 2008, Public Law 110-425 (H.R. 6353).

\bibitem{DowHC16}
Dowell, D., Haegerich, T.M., Chou, R.
CDC guideline for prescribing opioids for chronic pain--United States, 2016.
{\em MMWR Recomm. Rep.} 65, 1--49 (2016).

\bibitem{BefA15a} 
Befekadu, G.K., Antsaklis, P.J.
{\em On the asymptotic estimates for exit probabilities and minimum exit rates of diffusion processes pertaining to a chain of distributed control systems},
SIAM J. Contr. Optim. 53, 2297--2318 (2015).

\bibitem{DelM10}
Delarue, F., Menozzi, S. 
{\em Density estimates for a random noise propagating through a chain of differential equations}.
J. Funct. Anal. 259, 1577--1630 (2010).

\bibitem{FleSh85}
Fleming, W.H., Sheu, S.J.
{\em Stochastic variational formula for fundamental solutions of parabolic pde}.
Appl. Math. Optim. 13, 193--204 (1985).

\bibitem{She91}
Sheu, S.J.  
{\em Some estimates of the transition density of a nondegenerate diffusion Markov process}.
Ann. Probab. 19,  538--561 (1991).

\bibitem{Hor67}
H\"{o}rmander, L.
{\em Hypoelliptic second order differential operators},
Acta Math. 119, 147--171 (1967).

\bibitem{Ell73}
Elliott, D.L. 
{\em Diffusions on manifolds arising from controllable systems},
in Geometric Methods in System Theory, D.~Q. Mayne and R.W.~Brockett, eds. Reidel Publ. Co., Dordrecht, Holland, pp.~285--294 (1973).

\bibitem{FreWe84}
Freidlin, M.I., Wentzell, A.D. 
{\em Random perturbations of dynamical systems}.
Springer, Berlin (1984).

\bibitem{VenFre70}
Ventcel, A.D., Freidlin, M.I. 
{\em On small random perturbations of dynamical systems}.
Russian Math. Surv. 25, 1--55 (1970).

\bibitem{Aro67}
Aronson, D.G.
{\em Bounds for the fundamental solution of a parabolic equation}.
Bull. Amer. Math. Soc. 73, 890--896 (1967).

\bibitem{Roy88}
Royden, H.L. 
Real analysis.
Prentice Hall, Englewood Cliffs, NJ (1988).

\end{thebibliography}
\end{document}